\newcommand{\Z}{{\mathbb Z}}

\newcommand{\R}{{\mathbb R}}
\newcommand{\D}{{\mathbb D}}
\newcommand{\mbar}{{\overline{M}}}

\newcommand{\twid}{\widetilde}
\newcommand{\iso}{\approx}

\def\co{\colon\thinspace}
\documentclass[12pt]{amsart}
\usepackage{amssymb}
\usepackage{amsthm}
\usepackage{amsmath}
\usepackage{amscd}
\usepackage{graphicx}
\usepackage[all]{xy}
\usepackage[dvips,colorlinks=true]{hyperref}
\hypersetup{%
    pdftitle={Cohomology Determinants of Compact 3-Manifolds},
    pdfauthor={Christopher Truman},
    pdfkeywords={pdf, latex, tex, ps2pdf, dvipdfm, pdflatex},
    bookmarksnumbered,
    pdfstartview={FitH},
    urlcolor=cyan,
}%

\newtheorem{theorem}{Theorem}[section]
\newtheorem{proposition}[theorem]{Proposition}
\newtheorem{lemma}[theorem]{Lemma}
\newtheorem*{claim}{Claim}
\DeclareMathOperator{\Det}{Det}
\DeclareMathOperator{\Tors}{Tors}
\DeclareMathOperator{\Hom}{Hom}
\DeclareMathOperator{\im}{im}

\DeclareMathOperator{\rank}{rank}

\begin{document}
\title[Cohomology Determinants]{Cohomology determinants of compact $3$--manifolds}
\author{Christopher Truman}
\address{Mathematics Department\\University of Maryland\\College Park, MD  20742 USA}
\email{cbtruman@math.umd.edu}
\urladdr{http://www.math.umd.edu/\textasciitilde cbtruman/}
\begin{abstract}
We give definitions of cohomology determinants for compact, connected, orientable $3$--manifolds.  We also give formulae relating cohomology determinants before and after gluing a solid torus along a torus boundary component.  Cohomology determinants are related to Turaev torsion, though the author hopes that they have other uses as well.
\end{abstract}
\date{\today}
\maketitle
\begin{section}{Introduction}
\label{section-Introduction}
Cohomology determinants are an invariant of compact, connected, orientable $3$-manifolds.  The author first encountered these invariants in \cite{Tur:3d}, when Turaev gave the definition for closed $3$--manifolds, and used cohomology determinants to obtain a leading order term of Turaev torsion.  In \cite{Truman:TorsCohom}, the author gives a definition for $3$--manifolds with boundary, and derives a similar relationship to Turaev torsion.  Here, we repeat the definitions, and give formulae relating the cohomology determinants before and after gluing a solid torus along a boundary component.  One can use these formulae, and gluing formulae for Turaev torsion from \cite{Tur:3d} Chapter~VII, to re-derive the results of \cite{Truman:TorsCohom} from the results of \cite{Tur:3d} Chapter~III, or vice-versa. 
\end{section}
\begin{section}{Integral Cohomology Determinants}
\label{section-IntCohom}
%
\begin{subsection}{Closed $3$--manifolds}
\label{subsection-det-closed}
We will simply state the relevant result from \cite{Tur:3d} Section~III.1; the proof is similar to the one below.  Let $R$ be a commutative ring with unit, and let $N$ be a free $R$--module of rank $n\geq 3$.  Let $S=S(N^*)$ be the graded symmetric algebra on $N^*=\Hom_R(N,R)$, with grading $S=\bigoplus\limits_{\ell \geq 0} S^\ell$.  Let $f\co N\times N\times N\longrightarrow R$ be an alternate trilinear form, and let $g\co N\times N\rightarrow N^*$ be induced by $f$.  Now some notation: for $\theta$ any $(n\times n)$ matrix, we will denote by $\theta(i;j)$ the $(n-1 \times n-1)$ matrix obtained by striking the $i^\text{th}$ row and $j^\text{th}$ column from $\theta$.  Also, $[a'/a] \in R^\times$ is used to denote the determinant of the change of basis matrix from $a$ to $a'$.
\begin{lemma}[Turaev]
Let $\{a_i\}_{i=1}^{n}$, $\{b_i\}_{i=1}^{n}$ be bases of $N$ with dual basis $\{a_i^*\}_{i=1}^{n}$ of $N^*$.  Let $\theta$ be the $(n\times n)$ matrix over $S$ whose $i,j$ entry $\theta_{i,j}$ is given by $\theta_{i,j}=g(a_i,b_j)$.  Then there is a unique $d=d(f,a,b)\in S^{n-3}$ such that for any $1\leq i,j\leq n$,
\[ \det(\theta(i;j)) = (-1)^{i+j}a_i^*b_j^* d.\]
For any bases $a,a',b,b'$, we have
\[ d(f,a',b') = [a'/a][b'/b]d(f,a,b).\]
\end{lemma}

If $R=\Z$, then $d(f,a,a)$ is independent of the basis $a$, and will be denoted $\Det(f)$.  Now if $M$ is a closed, connected, oriented $3$--manifold with $b_1(M)\geq 3$, then $H^1(M)$ will satisfy the conditions for the module $N$ above, and the form $f_M\co H^1(M)\times H^1(M)\times H^1(M)\longrightarrow \Z$ given by $f_M(x,y,z)=\langle x\cup y\cup z , [M]\rangle$ is an alternate trilinear form.  Thus $f_M$ has a determinant $\Det(f_M)$.  This is what we will call the cohomology determinant.
\end{subsection}
\begin{subsection}{$3$--manifolds with nonvoid boundary}
\label{subsection-det-nonclosed}
Let $R$ be a commutative ring with unit, and let
$K,L$ be finitely generated free $R$ modules of rank $n$ and $n-1$
respectively, where $n \geq 2$.  As above, let $S =
S(K^*)$ where $K^* = \Hom_R(K,R)$.  Note if $\{a_i^*\}_{i=1}^{n}$ is
the basis of $K^*$ dual to the basis $\{a_i\}_{i=1}^{n}$ of $K$ then
$S=R[a_1^*,\dots,a_{n}^*]$, the polynomial ring on
$a_1^*,\dots,a_{n}^*$, and the grading of $S$ corresponds to the usual
grading of a polynomial ring.  Let $f\co L\times
K\times K \longrightarrow R$ be an $R$--module homomorphism which is
skew-symmetric in the two copies of $K$; i.e.  for
all $y,z \in K, x\in L, f(x,y,z)=-f(x,z,y)$.  Let $g$ denote the
associated homomorphism $L\times K\longrightarrow K^*$ given by
$(g(x,y))(z) = f(x,y,z)$.  
Again we set some notation: for a
matrix $A$, we will let $A(i)$ denote the matrix obtained by striking
out the $i^{\text{th}}$ column.
\begin{lemma}
\label{lemma-IntCohomDet}
Let
$\{a_i\}_{i=1}^n,\{b_j\}_{j=1}^{n-1}$ be bases for $K,L$ respectively,
and let $\{a_i^*\}$ 
be the basis of $K^*$ dual to the basis $\{a_i\}$ as above.  Let $\theta$ denote the $(n-1 \times n)$ matrix over $S$ whose
$i,j^{\text{th}}$ entry $\theta_{i,j}$ is obtained by
$\theta_{i,j} = g(b_i,a_j)$.  Then there is a unique $d = d(f,a,b) \in
S^{n-2}$ such that for any $1 \leq i \leq n$,
\begin{equation}
\det\ \theta(i) = (-1)^{i}a_i^*d.
\label{equation-IntCohomDetDef}
\end{equation}
For any other bases $a',b'$ of $K,L$ respectively, we have 
\begin{equation}
d(f,a',b') = [a'/a][b'/b]d(f,a,b).
\label{equation-IntCohomCob}
\end{equation}
\end{lemma}
\begin{proof}
Let $\beta$ denote the $(n-1\times n)$ matrix with $\beta_{i,j} =
g(b_i,a_j)a_j^*$.  The sum of the columns of $\beta$ is zero; indeed, for any
$i$, the $i^{\text{th}}$ entry (of the column vector obtained by summing the columns of $\beta$) is given by:
\begin{equation*}
\sum\limits_{j=1}^n \beta_{i,j} = \sum\limits_{j=1}^n g(b_i,a_j)a_j^* =
\sum\limits_{j,k=1}^n f(b_i,a_j,a_k)a_j^*a_k^* = 0.
\end{equation*}
The last equality follows since the $f$ term is anti-symmetric in
$j,k$ and the $a$ term is symmetric.  We now claim
$(-1)^{i}\det\beta(i)$ is independent of $i$.  
\begin{claim}
Let $Z$ be
a $(n-1\times n)$ matrix with columns $c_i$ for $1 \leq i \leq n$ such that
$\sum\limits_{i=1}^n c_i = 0$.  Then $(-1)^i \det(Z(i))$ for $1 \leq i
\leq n$ is independent of $i$. 
\end{claim}
The proof of this claim is as follows: think of $\det$ as a function on
the columns; 
$\det(Z(i))=\det(c_1 , c_2 , \dots , c_{i-1} ,
c_{i+1}, \dots, c_n)$.  Let $k\neq i$, then
$c_k = -\sum\limits_{p\neq k} c_p$, hence 
\begin{eqnarray*}
\det(Z(i))& = & \det(c_1,c_2, \dots , c_{k-1} , -\sum\limits_{p\neq
k} c_p, c_{k+1}, \dots ,  c_{i-1} ,
c_{i+1}, \dots , c_n)\\
& = & \sum\limits_{p\neq k} \det(c_1,c_2, \dots , c_{k-1} , -c_p,
c_{k+1}, \dots ,  c_{i-1} , c_{i+1}, \dots , c_n)\\
& = & \det(c_1,c_2, \dots , c_{k-1} , -c_i,
c_{k+1}, \dots ,  c_{i-1} , c_{i+1}, \dots , c_n).
\end{eqnarray*}
Here for notational convenience we have
assumed $k < i$, but it clearly makes no difference.  The last
equality holds because in each term but the $p=i$ term, we will have
two columns appearing twice.  Now to move $-c_i$ to the $i^\text{th}$
column, we will have to do $i-k-1$ column swaps.  Doing 
the column swaps and accounting for the negative sign of $c_i$, we get
$\det(Z(i))=(-1)^{i-k}\det(Z(k))$, which completes the proof of our
claim. 

This means
$(-1)^{i}\det\beta(i)$ is independent of $i$.  Now let $t_i = \det
\theta(i) \in S^{n-1}$.  It is clear that
\begin{equation*}
\det (\beta(i)) = t_i \prod_{k \neq i}a_k^*.
\end{equation*}
Then for any $i,p \leq n$, we have
\begin{eqnarray*}
(-1)^{i}t_ia_p^*\prod_{k=1}^n a_k^* & = & (-1)^{i}\det \beta(i)
a_p^*a_i^* \\
& = & (-1)^p\det \beta(p) a_i^*a_p^*\\ 
& = & (-1)^pt_pa_i^*\prod_{k=1}^n a_k^*.
\end{eqnarray*}
Now since the annihilators of $a_k^*$ in $S$ are zero, we must have 
\begin{equation*}
(-1)^{i}t_ia_p^* = (-1)^pt_pa_i^*.
\end{equation*}
This means that $a_i^*$ divides $t_ia_p^*$ for all $p$, hence $a_i^*$
divides $t_i$.  Define $s_i$ by $t_i = s_ia_i^*$.  Note 
\begin{equation*}
(-1)^{i}s_ia_i^*a_p^* = (-1)^{i}t_ia_p^* =
(-1)^pt_pa_i^*=(-1)^ps_pa_p^*a_i^*.
\end{equation*}
This means $(-1)^{i}s_i$ is independent of $i$.  Let $d =
(-1)^{i}s_i$.  By definition, 
\begin{equation*}
(-1)^{i}\det \theta(i) = (-1)^{i}t_i = (-1)^{i}s_ia_i^* = a_i^*d.
\end{equation*}
This proves \eqref{equation-IntCohomDetDef}.

Now to prove the change of basis formula, note we do not have to change both bases simultaneously, but can instead first obtain the formula for $d(f,a',b)$ in terms of $d(f,a,b)$, and then do the same for $b'$ and $b$.  So let $\{a_i'\}$ be another basis for $K$.  We
show $d(f,a',b) = [a'/a]d(f,a,b)$.  Let $S_i$ be the $(n\times n-1)$ matrix obtained by
inserting a row of zeroes into the $(n-1\times n-1)$ identity matrix
as the $i^\text{th}$ row.  Then one may easily see for any $(n-1\times
n)$ matrix $A$, the matrix $A(i)$ (obtained by striking out the
$i^\text{th}$ column) can also be obtained as $A(i)=AS_i$.  Let
$S_i^+$ denote the $(n\times n)$ matrix obtained by appending a column
vector with a $1$ in the $i^\text{th}$ entry and zeroes otherwise on
to the right of $S_i$, and let $A_+^i$ denote the $(n\times n)$ matrix
obtained by appending a row vector with a $1$ in the $i^\text{th}$
entry and zeroes otherwise on to the bottom of $A$.  Note
\begin{align*}
\det(S_i^+) &= (-1)^{n+i}\\
\det(A_+^i) &= (-1)^{n+i}\det(A(i))\\
&\text{hence}\\
\det(A_+^iS_i^+) &= \det(AS_i) = \det(A(i)).
\end{align*}
Now let $(a'/a)$
denote the usual change of basis matrix so that $a'_i =
\sum\limits_{j=1}^n (a'/a)_{i,j} a_j$.   Now $\theta_{i,j} =
g(b_i,a_j)$, so let 
\begin{align*}
\theta'_{i,j} &= g(b_i,a'_j)\\
 &= g(b_i,\sum\limits_{k=1}^n (a'/a)_{j,k} a_k)\\
 &= \sum\limits_{k=1}^n g(b_i,a_k)(a'/a)_{j,k}.
\end{align*}
Thus $\theta' = \theta\cdot (a'/a)^\mathsf{T}$.  Now 
\begin{align*}
\det\left(\theta_+^i (a'/a)^\mathsf{T} S_i^+\right) &=
\det\left(\theta_+^i\right)\det\left((a'/a)^\mathsf{T}\right)\det\left(S_i^+\right) \\
&=\det\left((a'/a)^\mathsf{T}\right) \det\left(\theta_+^i\right)\det\left(S_i^+\right) \\
&=\det\left(a'/a\right) \det\left(\theta_+^i\cdot S_i^+\right) \\
&= [a'/a]\det(\theta(i)) \\
&=[a'/a](-1)^i a_i^* d(f,a,b).
\end{align*}
Now we will compute the same thing in a much longer way to complete our proof.  Let
$e_i$ denote the row vector with a $1$ in the $i^\text{th}$ position and zeroes otherwise, i.e. the $i^\text{th}$ basis vector of $a$ as expressed in the $a$--basis, and let $r_i$ denote
the $i^\text{th}$ row of $(a'/a)^\mathsf{T}$ and $c_i$ denote the
$i^\text{th}$ column.  Then 
{\small 
\begin{align*}
\det(\theta_+^i (a'/a)^\mathsf{T} S_i^+) &= \det\left[\left(\begin{smallmatrix}
\theta \\ e_i\end{smallmatrix}\right) (a'/a)^\mathsf{T} \left(\begin{smallmatrix}
S_i & e_i^\mathsf{T} \end{smallmatrix}\right)\right] \\
&= \det\left[\left(\begin{smallmatrix}
\theta \\ e_i\end{smallmatrix}\right) \left(\begin{smallmatrix}
(a'/a)^\mathsf{T}(i) & c_i\end{smallmatrix}\right)\right] \\
&= \det\left(\begin{smallmatrix} \theta (a'/a)^\mathsf{T}(i) & \theta c_i
\\ (a'/a)^\mathsf{T}(i)_i & (a'/a)^\mathsf{T}_{i,i} \end{smallmatrix}\right) \\
&= (-1)^{n-i}\det\left(\begin{smallmatrix}\theta (a'/a)^\mathsf{T} \\
r_i\end{smallmatrix}\right) \\
&= (-1)^{n-i}\det\left(\begin{smallmatrix}\theta' \\
r_i\end{smallmatrix}\right) \\
&= (-1)^{n-i}\sum\limits_{k=1}^n
(-1)^{n+k}(a'/a)^\mathsf{T}_{i,k}\det(\theta'(k)) \\
&= (-1)^{n-i}\sum\limits_{k=1}^n (-1)^{n+k}(a'/a)^\mathsf{T}_{i,k} (-1)^k
(a'_k)^* d(f,a',b) \\
&= (-1)^{n-i}\sum\limits_{k=1}^n (-1)^{n+k}(a^*/(a')^*)_{i,k} (-1)^k
(a'_k)^* d(f,a',b) \\
&= (-1)^i d(f,a',b) a_i^*.
\end{align*}
}
So $d(f,a',b)-[a'/a]d(f,a,b)$ annihilates $a_i^*$ for each $i$, hence
is zero.

The computation for a $b$ change of basis is easier.  Let $b'$ be
another basis for $L$ and let $(b'/b)$ denote the $b$ to $b'$ change
of basis matrix.  Let $\theta'$ denote the matrix $g(b_i',a_j)$, then
$\theta' = (b'/b)\theta$.  So $\theta' S_i = (b'/b)\theta S_i$, hence
$\det(\theta'(i)) = [b'/b]\det(\theta(i))$.
This proves $d(f,a,b') = [b'/b]d(f,a,b)$, and completes the proof of
\eqref{equation-IntCohomCob}. 
\end{proof}

In the case, $R=\Z$,  our determinant depends on the basis only by
its sign.  In this case, we can refine the determinant by a choice of
orientation of the $\R$--vector space $(K\oplus L)\otimes \R$.  Let
$\omega$ be such a choice of orientation.  Then define $\Det_\omega(f)
= \det(f,a,b)$ where $a,b$ are bases of $K,L$ respectively such that
the induced basis of $(K\oplus L)\otimes \R$ given by
$\{a_1\otimes~1,a_2\otimes~1,\dots,a_n\otimes~1,b_1\otimes~1,b_2\otimes~1,\dots,b_{n-1}\otimes~1\}$
is positively oriented with respect to 
$\omega$.  Then $\Det_\omega(f)$ is well defined, and for any bases
$a',b'$, we have $\det(f,a',b') = \pm \Det_\omega(f)$ where the $\pm$
is chosen depending on whether $a',b'$ induces a positively or
negatively oriented basis of $(K\oplus L)\otimes \R$ with respect to
$\omega$.  

If $M$ is a compact, connected, oriented $3$--manifold with $\partial M\neq \varnothing$, $\chi(M)=0$, and $b_1(M)\geq 2$, then $H^1(M)$ and $H^1(M,\partial M)$ satisfy the conditions for $L,K$ respectively, and the map $f_M\co H^1(M,\partial M)\times H^1(M)\times H^1(M)\longrightarrow \Z$ given by $f_M(x,y,z)=\langle x\cup y\cup z , [M]\rangle$ has a determinant.  Also, a choice
of homology orientation of $M$ will determine an orientation for $(K\oplus
L)\otimes \R$.  To see why, let $\omega$ be a homology orientation
for $M$.  Consider $\{a_1^*,\dots,a_n^*\}$ a basis for
$H^1(M;\R)$ dual to a basis $\{a_1,\dots,a_n\}$ of $H_1(M,\R)$, and 
$\{b_1^*,\dots,b_{n-1}^*\}$ a basis of $H^1(M,\partial M;\R)$.  We will
say what it means for $\{a_1^*,\dots,a_n^*,b_1^*,\dots,b_{n-1}^*\}$ to be
a positively oriented basis for $H^1(M;\R)\oplus H^1(M,\partial M;\R)$,
and this will define our orientation.  Let
$[M]$ denote the fundamental class of $M$ determined by the
orientation of $M$ (not the homology orientation).  Then we will
define an orientation of $H^1(M;\R)\oplus H^1(M,\partial M;\R)$ by saying
that $\{a_1^*,\dots,a_n^*,b_1^*,\dots,b_{n-1}^*\}$ is a positively
oriented basis if and only if $\{[\text{pt}],a_1,\dots,a_n,b_1^*\cap
[M],\dots,b_{n-1}^*\cap [M]\}$ is a positively oriented basis for
$H_*(M;\R)$ with respect to $\omega$.  We will denote the sign
refined determinant with respect to this orientation by $\Det_{\omega}(f_M)$.  This is what we will call the cohomology determinant for a $3$--manifold with boundary.
\end{subsection}
\end{section}
\begin{section}{Gluing Formulae}
\label{section-Gluing}
We now give formulae for cohomology determinants for gluing solid tori to $3$--manifolds with nonempty boundary homeomorphic to a disjoint union of tori.  The reason we are only interested when the boundary consists entirely of tori is simple:  if $M$ is a compact orientable $3$--manifold with $\chi(M)=0$ and $\partial M$ containing a component not homeomorphic to a torus, then $\partial M$ contains a component homeomorphic to a sphere $S^2$.  This implies that the determinant for $M$ is zero.
\begin{subsection}{Gluing Homology Orientations}
\label{subsection-GluingHomologyOrientations}
To have a well-defined sign in our gluing formula, we must first review how to glue homology orientations.  The following is based on \cite{Tur:3d} Chapter~V, with some changes in notation.  Also, for simplicity, we will consider the solid tori being glued one-at-a-time, i.e. we will only give the definition for gluing one solid torus, and will consider the definition for gluing multiple solid tori to be given inductively; we can do this from \cite{Tur:3d} Lemma~V.2.3.  First, we define a {\it directed solid torus} as a solid torus $Z = \D^2\times S^1$ (where $\D^2$ is the standard 2-disk) with a distinguished generator of $H_1(Z)\iso\Z$, i.e. an orientation of the core $S^1$.  Now if $M$ is a compact connected 3-manifold with boundary consisting of tori and one boundary component $T$ picked out, then we can consider $\mbar=M\cup_T Z$ (under some choice of homeomorphism $T\rightarrow \partial Z$).  We can consider $Z$ to be homology oriented by setting $\omega_Z$ to be the orientation of $([\text{pt}],d)$ where $d$ is the distinguished generator of $H_1(Z)$ (to be precise, we should note that we are extending scalars from $\Z$ to $\R$).  This provides $H_*(Z,\partial Z;\R)$ with an orientation via Poincar\'{e} duality by saying $a\in H_2(Z,\partial Z)$ and $b\in H_3(Z,\partial Z)$ give a positively oriented basis $(a,b)$ of $H_*(Z,\partial Z)$ if and only if $(b^*\cap[Z],a^*\cap[Z])$ is a positively oriented basis of $H_*(Z)$ where $[Z]$ is either orientation class of $Z$.  It is clear that the resulting homology orientation of $H_*(Z,\partial Z;\R)$ does not depend on the (arbitrarily) chosen orientation of $Z$, but only depends on the distinguished direction of $Z$ (i.e. the distinguished generator of $H_1(Z)$).  This then provides $H_*(\mbar,M)$ with an orientation via excision; denote this orientation $\omega_{(\mbar,M)}$.  Then we may define $\twid{\omega},$ an orientation of $H_*(\mbar),$ from a given homology orientation $\omega$ of $M$ and our earlier constructed $\omega_{(\mbar,M)}$.  We define the orientation $\twid{\omega}$ of $H_*(\mbar)$ by requiring that the torsion of the homology exact sequence with $\R$ coefficients of the pair $(\mbar,M)$ have a positive sign (this exact sequence is an acylic complex).  Then we define the homology orientation of $\mbar$ induced from $\omega$, $\omega^\mbar$, as
\begin{equation}
\label{equation-hom-orien-glue-def}
\omega^\mbar = (-1)^{b_3(\mbar) + (b_1(M) + 1)(b_1(\mbar)+1)}\twid{\omega}.
\end{equation}
The sign in the equation is needed to guarantee that if we use this definition multiple times to glue on several directed solid tori, that the end result is independent of the order in which we perform our gluing, see \cite{Tur:3d} Lemma~V.2.3. 
\end{subsection}
\begin{subsection}{The Gluing Formulae}
\label{subsection-actual-gluing-formula}
\begin{theorem}
\label{theorem-dets-glue}
Let $M$ be a compact, connected, oriented 3-dimensional manifold with nonempty boundary consisting of tori and homology orientation $\omega$.  Let $\mbar$ be obtained by gluing a directed solid torus $Z$ along one boundary component $T$ of $M$, and let $\ell$ denote the image in $S(H_1(\mbar)/\Tors(H_1(\mbar)))$ of the distinguished generator of $H_1(Z)$.  If $\mbar$ is closed, assume $b_1(\mbar)\geq 3$, and if not assume $b_1(\mbar)\geq 2$.  
\begin{enumerate}
\item\label{item-det-glue-1} If $\partial M\neq T$ and the image of $H_1(T)$ in $H_1(M)$ is not rank 2, then \[\Det_\omega(f_M)=0.\]
\item\label{item-det-glue-2} If $\partial M=T$ and $b_1(\mbar)\neq b_1(M)$ then \[(\imath_M)_*(\Det_\omega(f_M))=0.\]
\item\label{item-det-glue-3} If $\partial M\neq T$ and the image of $H_1(T)$ in $H_1(M)$ is of rank 2, then 
\begin{equation*}
|\Tors(H_1(M))|(\imath_M)_*(\Det_\omega(f_M)) = |\Tors(H_1(\mbar))|\cdot\ell\cdot\Det_{\omega^\mbar}(f_\mbar).
\end{equation*}
\item\label{item-det-glue-4} If $\partial M=T$ and $b_1(\mbar)=b_1(M)$ then let $s_0$ denote the sign of the orientation $\omega^\mbar$ with respect to the natural homology orientation of $\mbar$ induced by an orientation.  Then
\begin{equation*}
|\Tors(H_1(M))|(\imath_M)_*(\Det_\omega(f_M)) = s_0 |\Tors(H_1(\mbar))|\cdot\ell\cdot\Det(f_\mbar).
\end{equation*}
\end{enumerate}
\end{theorem}
\end{subsection}
%
\begin{subsection}{Preliminary Remarks}
\label{subsection-prelim-remarks}
Let $M$ be a compact, connected, oriented $3$--manifold with $\partial M=\coprod\limits_{i} T_i$ where the index $i$ runs over some nonempty finite set, and each $T_i$ is a torus.  We will also denote $T=T_1$ and $R=\coprod\limits_{i > 1} T_i$ so that $\partial M=T\coprod R$ (note if $\partial M$ has one component $T$, then $R=\varnothing$).  We will be gluing a solid torus along the boundary component $T$ and will use $\mbar$ to denote the result, i.e. $\mbar=M\bigcup\limits_{T} Z$ for a solid torus $Z$ (the actual homeomorphism of $T$ to $\partial(\D^2\times S^1)$ will of course matter in the actual construction of $\mbar$, but we will not include it in our notation for simplicity).  We will also assume that $M$, $\mbar$ are given consistent orientations.  Since there is a difference in definition of the determinant for $\mbar$ closed, we will study the cases $R\neq\varnothing$ and $R=\varnothing$ separately.  Here let us also set some notation for the rest of the section.  We will often let $\lambda,\mu$ be a basis of $H_1(T)$ such that $\mu$ is the curve along which we will glue the meridian of our solid torus and $\lambda$ is parallel to the distinguished generator of $H_1(Z)$.  In other words, $\mu$ is killed in $H_1(\mbar)$, and $\lambda$ maps to $h\in H_1(\mbar)$.  The assumptions $b_1(\mbar)\geq 2$ if $\partial\mbar\neq\varnothing$ and $b_1(\mbar)\geq 3$ if $\partial\mbar=\varnothing$, will guarantee the appropriate ranges for $b_1(M)$ so that we will have well defined determinants for both $M$ and $\mbar$.

Whether $\mbar$ is closed or not, we must analyze mappings in cohomology; there is an obvious and natural map $H^1(\mbar)\rightarrow H^1(M)$ induced by the inclusion $M\hookrightarrow\mbar$.  However, $\partial M$ does not map to $\partial\mbar$ under the inclusion, so it does not induce a map from $H^1(\mbar,\partial\mbar)$ to $H^1(M,\partial M)$.  This means we will require a way to work around this unfortunate detail.  

Before we do so, however, we will give some results that we will be using throughout the section.  First, note by excision:
\begin{equation}
\label{equation-H^i(mbar,M)-excision}
H^i(\mbar,M)\iso H^i(Z,\partial Z)\iso \left\{\begin{array}{ll}\Z & \text{\ if\ } i=2,3\\0 & \text{\ otherwise.\ }\end{array}\right.
\end{equation}
Combining \eqref{equation-H^i(mbar,M)-excision} with the cohomology exact sequence of the pair $(\mbar , M)$
\begin{equation}
\label{equation-(mbar,M)-seq}
H^1(\mbar,M)\rightarrow H^1(\mbar)\rightarrow H^1(M)\rightarrow H^2(\mbar,M)
\end{equation}
we see that the cokernel of $H^1(\mbar)\rightarrow H^1(M)$ is rank $0$ or rank $1$, and the kernel is 0.  This means $b_1(\mbar)$ can either be $b_1(M)$ or $b_1(M)-1$.  Intuitively, the two cases correspond to either killing a finite order element or an infinite order element when we glue the solid torus along $T$.

We will also need to know something about how Poincar\'{e} duality compares before and after gluing.  Intuitively, one would expect that ``away from $T$'' (whatever that means), duality should be largely unchanged.  We now precisely state this intuitive idea.  To set some convenient notation, we will use $\imath_M$ to denote the inclusion $M\hookrightarrow \mbar$, $\imath_R$ to denote the inclusion $R\hookrightarrow\partial M$, and finally $\imath_{\partial\mbar}$ to denote the inclusion $\partial\mbar\rightarrow\partial M$ (by itself, this is the same as $\imath_R$, but we will use the notation $\imath_{\partial\mbar}$ when we want to look at induced maps for the triple $(\mbar,\partial M,\partial\mbar)$ and $\imath_R$ to look at induced maps for the triple $(M,\partial M,R)$).  Note that the map induced on cohomology by $\imath_M$ maps $H^*(\mbar,\partial\mbar)$ to $H^*(M,R)$.
\begin{proposition}
\label{proposition-takingCareOf[M],[mbar]}
For $w\in H^1(\mbar,\partial\mbar)$, if there is a $w'\in H^1(M,\partial M)$ such that $(\imath_M)^*(w) = (\imath_R)^*(w')\in H^1(M,R)$, then 
\[
w\cap [\mbar] =(\imath_M)_*(w'\cap [M]).
\]

\end{proposition}
Proposition~\ref{proposition-takingCareOf[M],[mbar]} will allow us to work around the fact that the inclusion $\imath_M\co M\hookrightarrow\mbar$ does not induce a map $\partial M\rightarrow \partial\mbar$, hence does not induce a map $H_3(M,\partial M)\rightarrow H_3(\mbar,\partial\mbar)$.  In particular, the inclusion does not induce anything so convenient as the map $[M]\mapsto [\mbar]$ of $H_3(M,\partial M)\rightarrow H_3(\mbar,\partial\mbar)$.  Furthermore, $\imath_M$ does not induce a nice map $H^1(\mbar,\partial\mbar)\rightarrow H^1(M,\partial M)$, so this Proposition helps us work around that as well.
\begin{proof}
First, write the commutative ladder induced by the inclusion $M\hookrightarrow \mbar$ with rows given by the cohomology sequences of the triples $(M,\partial M,R)$ and $(\mbar,\partial M,\partial\mbar)$ (note $\partial\mbar=R$ and could be empty):
\begin{equation}
\label{equation-(mbar,M,R)-diagram}
\text{{\footnotesize$
\begin{CD}
   @.    H^2(\mbar,M) @=    H^2(\mbar,M)   @.     @.   \\
    @.  @AAA               @AAA                  \\
   H^0(\partial M,R)     @>>>  H^1(M,\partial M)         @>>>  H^1(M,R) @>>> H^1(\partial M,R)  \\
      @|              @AAA              @AAA    @|\\
  H^0(\partial M,\partial\mbar)     @>>>  H^1(\mbar,\partial M) @>>>  H^1(\mbar,\partial\mbar) @>>> H^1(\partial M,\partial\mbar) \\
   @.   @AAA               @AAA                  @.\\
   @.    0            @.    0              @.       \\
\end{CD}
$}}
\end{equation}
A simple diagram chase, assuming a suitable $w'$ exists, shows that there is a $\twid{w}\in H^1(\mbar,\partial M)$ mapping to $w\in H^1(\mbar,\partial\mbar)$ and $w'\in H^1(M,\partial M)$. Now note by Alexander duality, $H_3(\mbar,\partial M)$ is free of rank $2$, and we have the following diagram with any straight line exact:
\begin{equation}
\label{equation-H_3(mbar,bd M)-diagram}
{\footnotesize
\xymatrix{
H_2(\partial M,\partial\mbar) & & \ar[dl]_{(\imath_M)_*} H_3(M,\partial M) \ar[r] & H_2(\partial M)\ar@{=}[dl] \ar[r]  & H_2(M)\ar[dl]\\
& \ar[ul] \ar[dl] H_3(\mbar,\partial M) \ar[r] & H_2(\partial M) \ar[r] & H_2(\mbar) & \\
H_3(\mbar,M)& & \ar[ul]^{(\imath_{\partial\mbar})_*} H_3(\mbar,\partial\mbar)\ar[r] & H_2(\partial\mbar)\ar[ul]\ar[r] & H_2(\mbar)\ar@{=}[ul]\\
}}
\end{equation}
We see $H_3(\mbar,\partial M)$ is generated by the images of the orientation classes $[M]$ and $[\mbar]$, and the difference of those images maps to (plus or minus) the generator of $H_3(\mbar,M)$.  So we now perform some simple computations:
\begin{align*}
\twid{w}\cap (\imath_M)_*([M]) & = (\imath_M)_*\left( (\imath_M)^*(\twid{w})\right) \cap [M])\\
& = (\imath_M)_*(w'\cap[M]).\\
\twid{w}\cap (\imath_{\partial\mbar})^*([\mbar]) & = (\imath_{\partial\mbar})_*\left( (\imath_{\partial\mbar})^*(\twid{w})\right) \cap [\mbar])\\
& = (\imath_{\partial\mbar})_*(w\cap[\mbar])\\
& = w\cap [\mbar].
\end{align*}
The last equality follows since the map induced by $(\imath_{\partial\mbar})$ on $H_2(\mbar)$ is equality in diagram \eqref{equation-H_3(mbar,bd M)-diagram}.

So we want to compute the cap product of $\twid{w}$ with the difference of $(\imath_M)_*([M])$ and $(\imath_{\partial \mbar})_*([\mbar])$ and show that it is zero.  The chain complex $C_*(\mbar,\partial M)$ consists of the chain complex $C_*(M,\partial M)$ with an additional two-handle and an additional three-handle, and the difference we are interested in is the class in $H_3(\mbar,\partial M)$ of the additional three-handle.  To compute the cap product of $\twid{w}$ with this homology class, we evaluate $\twid{w}$ on a 1-front face, and this is the coefficient of the 2-back face.  But each 1-front face of our 3-handle lies on $T$, and $\twid{w}\in H^1(\mbar,\partial M)$ means $\twid{w}$ is zero when restricted to $\partial M$, in particular when restricted to $T$.
\end{proof}

Recall that our determinants lie in the symmetric algebras $S=S\left( (H^1(M))^*\right)$ and $\overline{S}=S\left( (H^1(\mbar))^*\right)$ (for $M,\mbar$ respectively), so here we briefly comment on $S,\overline{S}$ and the map $S\rightarrow \overline{S}$ induced by the inclusion $M\hookrightarrow\mbar$.  First, the map $H^1(\mbar)\rightarrow H^1(M)$ induced by inclusion induces a dual map $(H^1(M))^*\rightarrow (H^1(\mbar))^*$, and if we think of $(H^1(M))^*$ as simply $H_1(M)/\Tors(H_1(M))$ and $(H^1(\mbar))^*$ as simply $H_1(\mbar)/\Tors(H_1(\mbar))$, then the map $S\rightarrow \overline{S}$ is the map induced by $H_1(M)\rightarrow H_1(\mbar)$ (which maps $\Tors(H_1(M))\rightarrow \Tors(H_1(\mbar))$).  Now $H_1(M)\rightarrow H_1(\mbar)$ is onto (its cokernel is contained in $H_1(\mbar,M)=0$), and similarly $H^1(\mbar)\rightarrow H^1(M)$ is injective with free cokernel of rank $0$ or $1$.  If the cokernel is rank $0$, then $H^1(\mbar)\rightarrow H^1(M)$ and $H_1(M)/\Tors(H_1(M)) \rightarrow H_1(\mbar)/\Tors(H_1(\mbar))$ are isomorphisms, as is $S\rightarrow \overline{S}$.  If the cokernel is rank $1$, then we may choose a basis $\alpha_1,\dots,\alpha_{n-1}$ of $H^1(\mbar)$ and then construct a basis $a_1,\dots,a_n$ of $H^1(M)$ with $(\imath_M)^*(\alpha_i)=a_i$ for $1\leq i\leq n-1$, and $a_n$ having nonzero image in $H^2(\mbar,M)$.  Then the induced map $(H^1(M))^*\rightarrow (H^1(\mbar))^*$ is the map $a_i^*\mapsto\alpha_i^*$ for $1\leq i\leq n-1$, and $a_n^*\mapsto 0$ (and similarly $S\rightarrow\overline{S}$).  We will slightly abuse notation and denote the map $S\rightarrow \overline{S}$ by $(\imath_M)_*$.
\end{subsection}
\begin{subsection}{The Case $R\neq\varnothing$}
\label{subsection-R_nonempty}
In this case, we know that $\mbar$ is also a $3$--manifold with nonempty boundary, so we will use the determinant from \ref{section-IntCohom}.  First, a preliminary result involving rank counting.  
\begin{lemma}
\label{lemma-M,R-betti}
The following are all equal to zero:
\begin{equation}
\label{equation-M,R-betti-0,3=0}
b_0(M,T)=b_0(M,R)=b_3(M,T)=b_3(M,R)=0.
\end{equation}
The following are all equal:
\begin{equation}
\label{equation-M,R-betti-1,2}
b_1(M,T)=b_2(M,T)=b_1(M,R)=b_2(M,R).
\end{equation}
\end{lemma}
\begin{proof}
We first note $b_0(M,T)=0$ and $b_0(M,R)=0$ since $H_0(T)\rightarrow H_0(M)$ and $H_0(R)\rightarrow H_0(M)$ are both surjective, and then \[b_3(M,R)=b^3(M,R)=b_0(M,T)=0.\]  The first equation is by the universal coefficient theorem, the second is by Poincar\'{e} duality.  We similarly conclude $b_3(M,T)=0$.

Now $b_1(M,T)=b^2(M,R)=b_2(M,R)$ by the same reasoning as above, so it remains to show that $b_1(M,T)=b_2(M,T)$.  This follows from counting ranks in the exact sequence of the pair $(M,T)$ and noting that since $\chi(M)=\chi(T)=0$ and $b_0(M,T)=b_3(M,T)=0$, we must have $b_1(M,T)=b_2(M,T)$.
\end{proof}
Now we will look at (the first few terms of) the exact sequence of the triple $(M,\partial M,R)$ and the (reduced) exact sequence of the pair $(M,T)$ (both in cohomology):
\begin{equation}
\label{equation-(M,bdM,R)-seq}
0\rightarrow H^0(\partial M,R)\rightarrow H^1(M,\partial M)\rightarrow H^1(M,R)\rightarrow H^1(\partial M,R),
\end{equation}
\begin{equation}
\label{equation-(M,T)-seq}
0\rightarrow H^1(M,T)\rightarrow H^1(M)\rightarrow H^1(T)
.\end{equation}
Note that $H^1(\partial M,R)\iso H^1(T)$ and in fact we can form a commutative square with the last two terms each of \eqref{equation-(M,bdM,R)-seq} and \eqref{equation-(M,T)-seq}, where the right vertical arrow is an isomorphism:
\begin{equation}
\label{equation-M-R-T-square}
\begin{CD}
H^1(M,R) @>>> H^1(\partial M,R)\\
@VVV @VVV\\
H^1(M) @>>> H^1(T).
\end{CD}
\end{equation}
Since $H^1(T)\iso \Z^2$, the maximum rank of the image of each horizontal arrow is two, and by commutativity and the fact that the right vertical arrow is an isomorphism, the rank of the image of $H^1(M,R)$ in $H^1(\partial M,R)$, which we will denote by $s=\rank_\Z(\im(H^1(M,R)\rightarrow H^1(\partial M,R)))$, is less than or equal to the rank of the image of $H^1(M)$ in $H^1(T)$, which we will denote by $r=\rank_\Z(\im(H^1(M)\rightarrow H^1(T)))$ (i.e. $r\geq s$).  Now if $n=b_1(M)$ then $n-1=b_1(M,\partial M)$.  Note $H^0(\partial M,R)\iso\Z$ so counting ranks in \eqref{equation-(M,bdM,R)-seq} gives
\begin{equation}
\label{equation-b_1(M,R)-vs-rk}
b_1(M,R) = n-2+s.
\end{equation}  
\begin{subsubsection}{Case 1: $r=2$}
\label{subsubsection-case_1:r=2}
First, note that this can occur; for example the exterior of the Hopf link, where $T$ is either boundary component, has $H^1(M)\rightarrow H^1(T)$ an isomorphism.  So this case is not vacuous.

In this case, $b_1(M,T)=n-2$, so by \eqref{equation-b_1(M,R)-vs-rk} and Lemma~\ref{lemma-M,R-betti}, $s=0$.  Each group in both \eqref{equation-(M,bdM,R)-seq} and \eqref{equation-(M,T)-seq} is free, and $H^1(M,\partial M)\rightarrow H^1(M,R)$ is onto hence splits, so given a basis $\beta_1,\dots,\beta_{n-2}$ of $H^1(M,R)$, we may choose a basis $b_1,\dots ,b_{n-1}$ of $H^1(M,\partial M)$ such that $b_i\mapsto \beta_i$ for $1\leq i \leq n-2$ and $b_{n-1}$ is dual (under evaluation) to a path connecting $T$ to one of the components of $R$.  If we let $\imath_T\co T\hookrightarrow M$ be the inclusion, then $(\imath_T)_*([T]) = b_{n-1}\cap [M]$.

We now similarly compare $H^1(M)$ to $H^1(M,T)$.  Since $r=2$, for any basis $\alpha_1,\dots,\alpha_{n-1}$ of $H^1(M,T)$, we can choose a basis $a_1,\dots,a_n$ of $H^1(M)$ with $\alpha_i\mapsto a_i$ for $1\leq i\leq n-2$, and $a_{n-1},a_n$ mapping to linearly independent elements in $H^1(T)$.  Thus if we choose any basis $c_1,c_2$ of $H^1(T)$, then $\imath_T^*(a_{n-1}) = a_{1,1}c_1 + a_{1,2}c_2$ and $\imath_T^*(a_{n}) = a_{2,1}c_1 + a_{2,2}c_2$ where $A=\left( \begin{smallmatrix}a_{1,1} & a_{1,2}\\a_{2,1}&a_{2,2}\end{smallmatrix}\right)$ is an integral matrix with $\det(A) = D\neq 0$.

Now we look at the matrix over $S((H^1(M))^*)$ given by cup product as in Lemma~\ref{lemma-IntCohomDet}.  There will be a $n-2\times n-2$ square matrix in the upper left hand corner composed of the cup products of the $\alpha$'s and $\beta$'s, let us call this matrix $\mathfrak{M}$, and then the last two columns will be for cup products of the $b_i$ with $a_{n-1},a_n$ and the last row for $b_{n-1}$ cup the $a_j$.  Recall the matrix $\theta$ from Lemma~\ref{lemma-IntCohomDet}, $\theta_{i,j} = g(b_i,a_j)$; $\theta$ will have the form
\begin{equation}
\label{equation-r=2-R,T-cohom-matrix}
\left( \begin{array}{ccc}\text{{\Large{$\mathfrak{M}$}}} & v_{1} & v_{2}\\w & \pm Da_n^* & \mp D a_{n-1}^*\end{array}\right)
.\end{equation}
Above, $v_1,v_2$ are dimension $n-2$ column vectors and $w$ is a dimension $n-2$ row vector, and the signs are chosen depending on whether $c_1\cup c_2$ is dual, under evaluation, to $\pm [T]$.  Now since $(\imath_T)_*([T]) = b_{n-1}\cap [M]$, for any $u,v\in H^1(M)$, we can compute 
\[
\langle u\cup v\cup b_{n-1} , [M]\rangle = \langle u\cup v , b_{n-1}\cap [M]\rangle = \langle \imath_T^*(u)\cup\imath_T^*(v) , [T]\rangle.
\]
 This explains the $\pm D$ terms in the matrix, and also allows us to note that if $a_i\cup b_{n-1}\neq 0$ for some $i$, then there is some $v\in H^1(M)$ such that \[\langle a_i\cup v\cup b_{n-1},[M]\rangle \neq 0,\] so $\imath_T^*(a_i)\neq 0$.  This means that the row vector $w$ is equal to 0, since $\imath_T^*(\alpha_i)=0$ for $1\leq i \leq n-2$.  And now it is easy to compute the determinant, 
\begin{equation}
\label{equation-r=2-detM}
\det(\theta(n)) = (\pm D a_n^*) \det(\mathfrak{M}).
\end{equation}

Now we know $b_1(M,R)=n-2$, and we must have $b_1(\mbar,\partial\mbar)=n-2$ as well since we must have $b_1(\mbar)=b_1(M)-1$.  Geometrically, this means if each generator of $H_1(T)$ is infinite order in $H_1(M)$ (which corresponds to $r=2$), then by gluing a solid torus on $T$ we must kill an infinite order element.  Since $H^1(\mbar,\partial\mbar)$ injects into $H^1(M,R)$ with a free cokernel and they have the same rank, $H^1(\mbar,\partial\mbar)\rightarrow H^1(M,R)$ an isomorphism, 

We also look at the triple $(\mbar,M,T)$, using the following commutative diagram with exact rows and columns (we abuse notation and let $T$ denote the image of $T$ in $\mbar$):
\begin{equation}
\label{equation-(mbar,M,T)-diagram}
\begin{CD}
    @.    H^2(\mbar,M) @=    H^2(\mbar,M)   @.        \\
@.      @AAA               @AAA                  @.\\
  0 @>>>  H^1(M,T)     @>>>  H^1(M)         @>>>  H^1(T)  \\
@.      @AAA               @AAA                  @|\\
  0 @>>>  H^1(\mbar,T)     @>>>  H^1(\mbar) @>>>  H^1(T)  \\
@.      @AAA               @AAA                  @.\\
    @.    0            @.    0              @.        \\
\end{CD}
\end{equation}
This diagram is obtained by ``pulling apart'' (along the equalities) the braid diagram that gives rise to the exact sequence of the triple.  By commutativity, we note that the image of $H^1(\mbar)$ in $H^1(T)$ has rank $1$, and thus $H^1(\mbar,T)$ has rank $b_1(\mbar)-1 = n-2 = b_1(M,T)$ and the map $H^1(\mbar,T)\rightarrow H^1(M,T)$ is an injection with free cokernel of free $\Z$--modules of the same rank, hence is an isomorphism.  So we may choose a basis of $H^1(\mbar)$ by choosing a basis of $H^1(\mbar,T)$ and a preimage of the generator of the image of $H^1(\mbar)$ in $H^1(T)$, let us denote this element by $\alpha_{n-1}\in H^1(\mbar)$.

We now have chosen bases of $H^1(\mbar)$ and $H^1(\mbar,\partial\mbar)$ which are very similar to the bases of $H^1(M)$ and $H^1(M,\partial M)$, and the matrix we will want to study for the purposes of constructing the determinant, which we will call $\overline{\theta}$, will have the square $n-2\times n-2$ block in the upper left corner $(\imath_M)_*(\mathfrak{M})$ (this follows from Proposition~\ref{proposition-takingCareOf[M],[mbar]}).  This means 
\begin{equation}
\label{equation-det(overline(theta))}
\det(\overline{\theta}(n-1))=(\imath_M)_*(\det(\mathfrak{M})).
\end{equation}

Now recall our notation of $\lambda,\mu$ as the basis of $H_1(T)$ introduced in \ref{subsection-prelim-remarks}.  Then $\lambda^*,\mu^*$ is a basis of $H^1(T)$, and let $k\in\Z$ such that $\imath_T^*(\alpha_{n-1}) = k\lambda^*$ (we have no multiples of $\mu$ since $\langle \imath_T^*(\alpha_{n-1}),\mu\rangle = \langle \alpha_{n-1} , (\imath_T)_*(\mu)\rangle = 0$ since $\mu$ is killed in $\mbar$).  We can take $(\imath_M)^*(\alpha_{n-1})$ to be one of our generators in $H^1(M)$ with nonzero image in $H^1(T)$ by commutativity of \eqref{equation-(mbar,M,T)-diagram} and the fact that $(\imath_M)^*$ has free cokernel.  Choose any suitable $a_n$ for the final generator of $H_1(M)$, and let $m\in\Z$ such that $\langle a_n , (\imath_T)_*(\mu)\rangle = m$, so that the $D$ given in \eqref{equation-r=2-R,T-cohom-matrix} is simply $k\cdot m$ and note that since $\mu$ is killed in $H_1(\mbar)$, we are introducing new $m$--torsion to $H_1(\mbar)$, i.e. we have $|\Tors(H_1(\mbar))| = m\cdot|\Tors(H_1(M))|$.  Also, for simplicity, if the $D$ appearing in \eqref{equation-r=2-detM} is negative, we can change $a_n$ to $-a_n$ to force the sign of $D$ to be positive.

Now we are finally ready to compare the determinants of the forms $f_M$ and $f_\mbar$.  Let $a$ be the basis of $H^1(M)$ consisting of:
\[a_1=(\imath_M)^*(\alpha_1),\dots,a_{n-2}=(\imath_M)^*(\alpha_{n-2}),\] followed by $a_{n-1}=(\imath_M)^*(\alpha_{n-1})$ and then $a_n$.  Let $b$ be the basis of $H^1(M,\partial M)$ consisting of $b_1,\dots,b_{n-2}$ with $(\imath_R)^*(b_i)=(\imath_M)^*(\beta_i)$ (for $i\leq n-2$), followed by $b_{n-1}$.  Then with $\theta$ expressed in this basis, $\det(\theta(n)) = (-1)^n a_n^* d(f_M,a,b)$ by Lemma~\ref{lemma-IntCohomDet}.  But by \eqref{equation-r=2-detM}, we have (recalling we chose $a_n$ so that $D$ is positive)
\[
D a_n^* \det(\mathfrak{M}) = (-1)^n a_n^* d(f_M,a,b).
\]
This means
\begin{equation}
\label{equation-r=2,det(f_M)}
d(f_M,a,b) = (-1)^n D \det(\mathfrak{M})
.\end{equation}
Furthermore, by \eqref{equation-det(overline(theta))} and Lemma~\ref{lemma-IntCohomDet}, if we choose the basis $\alpha$ of $H^1(\mbar)$ to be $\alpha_1,\dots,\alpha_{n-2}$ followed by $\alpha_{n-1}$, and the basis $\beta$ of $H^1(\mbar,\partial\mbar)$ to be $\beta_1,\dots,\beta_{n-2}$, then
\begin{align}
\notag (\imath_M)_*(\det(\mathfrak{M})) &= \det(\overline{\theta}(n-1))\\
\label{equation-r=2,det(f_mbar)} & = (-1)^{n-1} \alpha_{n-1}^* d(f_\mbar,\alpha,\beta).
\end{align}
Now plugging \eqref{equation-r=2,det(f_mbar)} into \eqref{equation-r=2,det(f_M)} we obtain
\begin{equation*}
(\imath_M)_*(d(f_M,a,b))= -D \alpha_{n-1}^* d(f_\mbar,\alpha,\beta).
\end{equation*}
We have constructed this so that $\langle \alpha_{n-1} , \lambda\rangle = k$ meaning $k\alpha_{n-1}^*\mapsto \ell$ in the canonical isomorphism $(H^1(\mbar))^*\rightarrow H_1(\mbar)/\Tors(H_1(\mbar))$  where $\ell$ is the image of $\lambda$ in $H_1(\mbar)/\Tors(H_1(\mbar))$, so this means
\begin{equation}
\label{equation-det(f_M)-vs-det(f_mbar)}
(\imath_M)_*(d(f_M,a,b))= (-1) m\cdot \ell\cdot d(f_\mbar,\alpha,\beta).
\end{equation}
To complete the proof of Theorem~\ref{theorem-dets-glue} Item~\ref{item-det-glue-3}, we must see how the sign refined determinants work out using the induced homology orientation on $\mbar$.  To do so, we first take the sign of the torsion of the reduced long exact sequence of the pair $(\mbar,M)$.
We use the reduced sequence because we do not need the end of the sequence since it contributes nothing to the torsion.  Note the sign in \eqref{equation-hom-orien-glue-def} is trivial, so $\omega^\mbar$ is simply an orientation of $H_*(\mbar)$ making the torsion of the above sequence positive.  A simple calculation shows us that if $a,b,\alpha,\beta$ are bases as above and we use them to compute the sign of the torsion of the above sequence, we obtain a negative answer.  This means if $a$ and $b$ are bases such that $d(f_M,a,b)=\Det_\omega(f_M)$, then $\alpha$ and $\beta$ are bases such that $d(f_\mbar,\alpha,\beta)=-\Det_{\omega^\mbar}(f_\mbar)$.  This proves Theorem~\ref{theorem-dets-glue} Item~\ref{item-det-glue-3}.
\end{subsubsection}
\begin{subsubsection}{Case 2: $r=1$}
\label{subsubsection-case_2:r=1}
As in the earlier case, we first note that this case is not vacuous.  An example would be the exterior of two disjoint unlinked $S^1$'s embedded in $S^3$, with $T$ as either boundary component.  However, we will shortly see that the determinants in this case are as uninteresting as in our example.

We will first analyze the decompositions of $H^1(M,\partial M)$ and $H^1(M)$ with respect to $H^1(M,R), H^1(M,T), $ and $H^1(T)$ as before.  Note if $r=1$, then by \eqref{equation-(M,bdM,R)-seq} we know that $b_1(M,T)=n-1$ and hence $s=1$ as well.  So now choose $\alpha_1,\dots,\alpha_{n-1}$ a basis of $H^1(\mbar)$, and we can choose a basis of $H^1(M)$ with $\alpha_i\mapsto a_i$ for $1\leq i\leq n-1$ and $\imath_T^*(a_n)\neq 0$.

Now $H_1(M,\partial M)$ still has a free summand of rank $1$ generated by the dual of a path connecting $T$ to any component of $R$, but now $H^1(M,R)$ splits as the image of $H^1(M,\partial M)$ plus another free generator, which must map to $a_n$ under $H^1(M,R)\rightarrow H^1(M)$ since the cokernel of that map is free and everything else in $H^1(M,R)$ maps to zero in $H^1(T)$.  So choose bases $b_1,\dots,b_{n-1}$ of $H^1(M,\partial M)$ where $b_{n-1}$ is as above, dual to a path connecting $T$ to a component of $R$, and $b_i\mapsto \beta_i$ for $1\leq i\leq n-2$ where $\beta_1,\dots,\beta_{n-2},\gamma$ is a basis of $H^1(M,R)$ and $\gamma\mapsto a_n$ under $H^1(M,R)\rightarrow H^1(M)$.  Note we still have $(\imath_T)_*([T]) = b_{n-1}\cap [M]$.

In addition, we have, just as in our earlier case, for any $u,v\in H^1(M)$,
\[
\langle u\cup v\cup b_{n-1} , [M]\rangle = \langle \imath_T^*(u)\cup\imath_T^*(v) , [T]\rangle.
\]
Since $1=r=\rank_\Z(\im(\imath_T^*))$, we know $\imath_T^*(u),\imath_T^*(v)$ are both multiples of the same element in $H^1(T)$, so their cup product is zero.  This means that the last row of the matrix $\theta$ consists entirely of zeros, so $\det(\theta(i))=0$ for any $1\leq i\leq n$, proving Theorem~\ref{theorem-dets-glue} Item~\ref{item-det-glue-1}.
\end{subsubsection}
\begin{subsubsection}{Case 3: $r=0$}
\label{subsubsection:case_3:r=0}
Unlike the first two cases, this case is vacuous; it cannot occur since $r=0$ means $b_1(M,T)=b_1(M)=n$ hence $b_1(M,R)=n$, and then \eqref{equation-b_1(M,R)-vs-rk} gives $s=2$, contradicting our earlier claim that $r\geq s$.  Geometrically, this would correspond to the case that $H_1(T)\rightarrow H_1(M)$ has image entirely in $\Tors(H_1(M))$. One can verify that this cannot happen by letting $M'$ denote the result of gluing solid tori along each component of $R$ in any way one likes, thus $T=\partial M'$. Note the following commutative diagram (the cokernel of $H_1(M)\rightarrow H_1(M')$ is contained in $H_1(M',M)$, which is zero by the homology analogue of \eqref{equation-H^i(mbar,M)-excision}, and similarly for $H_1(M,T)\rightarrow H_1(M',T)$):
\[
\begin{CD}
     H_1(T)    @>>>  H_1(M)  @>>> H_1(M,T)  \\
      @|             @VVV          @VVV    \\
     H_1(T)    @>>>  H_1(M') @>>> H^1(M',T) \\
      @.             @VVV         @VVV   \\
               @.      0     @.       0
\end{CD}
\]
From this diagram, we note that the image of $H_1(T)$ cannot be rank $0$ in $H_1(M)$, because it is rank $1$ in $H_1(M')$.
\end{subsubsection}
\end{subsection}
\begin{subsection}{The Case $R=\varnothing$}
\label{subsection-R_empty}
In this case, we must compare the determinant from \ref{subsection-det-nonclosed} (when we are looking at $M$, before the gluing) to the determinant from \ref{subsection-det-closed} (i.e. from \cite{Tur:3d}~III.1).  First, we know $H^1(M,\partial M)\rightarrow H^1(M)$ is an injection with free cokernel, which must be of rank $1$ since $b_1(M,\partial M) = b_1(M)-1$; we will still use $n$ to denote $b_1(M)$.  Now we still have $H^i(\mbar,M)\iso H^i(\D^2\times S^1,S^1\times S^1)$, so we still have $b_1(\mbar)$ either equal to $n$ or $n-1$ depending on whether the element in $H_1(M)$ killed is finite or infinite order, and each of these cases can occur.  So let us examine both.  Also note that the reasoning from above (in \ref{subsubsection-case_2:r=1}) that led us to conclude that $\theta$ had a row consisting entirely of zeroes does not apply here, since that row corresponded to an element of $H^1(M,\partial M)$ connecting $T$ to another boundary component, and such a thing does not exist if $\partial M=T$ (in fact, the image of $[T]$ is zero in $H_2(M)$, so we will not be pulling back cohomology elements along the inclusion of $T$ into $M$ at all). 
\begin{subsubsection}{Case 1: $b_1(\mbar)=b_1(M)$}
\label{subsubsection-b_1(mbar)=n}
If $b_1(\mbar)=n$, then $(\imath_M)^*\co H^1(\mbar)\rightarrow H^1(M)$ is an isomorphism, and $H^1(M,\partial M)\rightarrow H^1(M)$ has kernel of rank $1$.  If we choose a basis $b_1,\dots,b_{n-1}$ of $H^1(M,\partial M)$, we can choose $a_n$ so that the images of the $b_i$, which we will call $a_i$, in $H^1(M)$ combined with $a_n$ forms a basis of $H^1(M)$, and we will let $\alpha_i\in H^1(\mbar)$ with $(\imath_M)^*(\alpha_i)=a_i$.  Then the matrix $(\imath_M)_*(\theta)$ will be all but the last row of the matrix $\overline{\theta}$ by Proposition~\ref{proposition-takingCareOf[M],[mbar]}, and $(\imath_M)_*(\det(\theta(n))) = \det(\overline{\theta}(n;n))$, and hence
\begin{align*}
(\imath_M)_*((-1)^n a_n^* d(f_M,a,b)) & = (\imath_M)_*(\det(\theta(n)))\\
& = \det(\overline{\theta}(n;n)) \\
& = \alpha_n^* \alpha_n^* d(f_\mbar,\alpha,\alpha).
\end{align*} 
Since $\alpha_n^*=(\imath_M)_*(a_n^*)$, the conclusion for determinants is that 
\[
(\imath_M)_*(d(f_M,a,b))= (-1)^n (\imath_M)_*(a_n^*)d(f_\mbar,a,a).
\]  
Now if we have chosen $\lambda,\mu$ as above, a basis of $H_1(T)$ so that $\mu$ is the basis element along which the meridian of our solid torus is glued, then $(\imath_T)_*(\mu)$ is finite order in $H_1(M)$ since gluing does not change the first Betti number, so let us say that $(\imath_T)_*(\mu)$ has order $k\in H_1(M)$; this means 
\[
\Tors(H_1(\mbar))\cdot k=\Tors(H_1(M)).
\]  
Then since $\Tors(H_1(M))\iso\Tors(H_1(M,\partial M))$, and $(\imath_T)_*(\mu)$ maps to zero in $H_1(M,\partial M)$, we must have a $k^{\text{th}}$ root of $(\imath_T)_*(\lambda)$ in $H_1(M)$, which we can choose $a_n$ to be dual to i.e. $\langle a_n , \lambda\rangle = k$.  Finally, if we let $\ell$ denote the image of $\lambda$ in $H_1(M)/\Tors(H_1(M))$, and since $(H^1(M))^*$ naturally isomorphic to $H_1(M)/\Tors(H_1(M))$, then we can write $k\alpha_n^*=\ell$, so multiplying through by $k$ we have
\begin{align*}
k\cdot (\imath_M)_*(d(f_M,a,b)) &= (-1)^n (k\alpha_n^*) d(f_\mbar,\alpha,\alpha)\\
& = (-1)^n \ell \cdot d(f_\mbar,\alpha,\alpha).
\end{align*}
To complete the proof of Theorem~\ref{theorem-dets-glue} Item~\ref{item-det-glue-4}, we must once again analyze signs.  First, the sign of $\omega^\mbar$ is equal to $(-1)^n$ times the sign of an orientation $\omega'$ of $H_*(\mbar)$ giving positive torsion of the exact sequence
\[ H_2(M)\rightarrow H_2(\mbar)\rightarrow H_2(M,\mbar)\rightarrow H_1(M)\rightarrow H_1(\mbar).\]
This time, we have truncated the sequence both on the left and right since the truncated parts did not contribute to the sign.  Another simple torsion calculation tells us that the sign $s_0$ in Theorem~\ref{theorem-dets-glue} Item~\ref{item-det-glue-4} is simply $(-1)^n$ times the sign of $d(f_M,a,b)$ with respect to $\Det_\omega(f_M)$.  This proves Theorem~\ref{theorem-dets-glue} Item~\ref{item-det-glue-4}.
\end{subsubsection}
\begin{subsubsection}{Case 2: $b_1(\mbar)=b_1(M)-1$}
\label{subsubsection-b_1(mbar)=n-1}
In this case, we may use the diagram \eqref{equation-(mbar,M,T)-diagram}, with $T=\partial M$, and we see that $H^1(\mbar,\partial M)\rightarrow H^1(M,\partial M)$ is an isomorphism, as is $H^1(\mbar,\partial M)\rightarrow H^1(\mbar)$.  So we may choose a basis $b_1,\dots,b_{n-1}$ of $H^1(M,\partial M)$ and additional element $a_n\in H^1(M)$ with the images of the $b_i$, which we call $a_i$, combined with $a_n$ is a basis of $H^1(M)$, and then $\overline{\theta} = (\imath_M)_*(\theta(n))$.  Now we can choose $\alpha_1,\dots,\alpha_{n-1}$ a basis of $H^1(\mbar)$ with $(\imath_M)^*(\alpha_i) = a_i$ for $1\leq i\leq n-1$.  Using this basis, we compute the determinant $\det(\theta(1)) = (-1)a_1^*d(f_M,a,b)$ by running down the $n^{\text{th}}$ column:
{\footnotesize
\begin{align*}
(\imath_M)_*(\det(\theta(1))) & = \sum\limits_{i < n} (-1)^{i+n}(\imath_M)_*(g(b_i,a_n))(\imath_M)_*(\det((\theta(1))(i,n)))\\
& = \sum\limits_{i < n} (-1)^{i+n}(\imath_M)_*(g(b_i,a_n))\det(\overline{\theta}(i;1))\\
& = \sum\limits_{i < n} (-1)^{i+n}(\imath_M)_*(g(b_i,a_n))(-1)^{i+1} \alpha_i^*\alpha_1^*d(f_\mbar,\alpha,\alpha)\\
& = (-1)^{n+1}\alpha_1^*d(f_\mbar,\alpha,\alpha)\sum\limits_{i < n} (\imath_M)_*(g(b_i,a_n))\alpha_i^*\\
& = (-1)^{n+1}(\imath_M)_*(a_1^*)d(f_\mbar,\alpha,\alpha)\sum\limits_{i < n} \sum\limits_{k=1}^{n-1}(\imath_M)_*(f_M(b_i,a_n,a_k)a_k^*)\alpha_i^*\\
& = (-1)^{n+1}(\imath_M)_*(a_1^*)d(f_\mbar,\alpha,\alpha)\sum\limits_{i < n} \sum\limits_{k=1}^{n-1}\langle b_i\cup a_n\cup a_k,[M]\rangle\alpha_k^*\alpha_i^*\\
& = (-1)^{n+1}(\imath_M)_*(a_1^*)d(f_\mbar,\alpha,\alpha)\sum\limits_{i=1}^{n-1} \sum\limits_{k=1}^{n-1}-\langle b_i\cup a_k, a_n\cap [M]\rangle\alpha_k^*\alpha_i^*\\
& = (-1)^{n+1}(\imath_M)_*(a_1^*)d(f_\mbar,\alpha,\alpha)\sum\limits_{i,k=1}^{n-1}-\langle b_i\cup b_k, a_n\cap [M]\rangle\alpha_k^*\alpha_i^*\\
& = 0.
\end{align*}}
The last equality is true since we are summing over $i,k$ and the $b_i\cup b_k$ is antisymmetric in $i,k$ and $\alpha_k^*\alpha_i^*$ is symmetric.  The line before that follows from noting that $a_k$ is the image of $b_k$ under $H^1(M,\partial M)\rightarrow H^1(M)$.This proves Theorem~\ref{theorem-dets-glue} Item~\ref{item-det-glue-2}, and in fact completes the proof of Theorem~\ref{theorem-dets-glue}.\qed
\end{subsubsection}
\end{subsection}
\end{section}
%
\begin{section}{Determinants and Massey Products}
\label{section-IntMass-Dets}
One can similarly construct a determinant using higher Massey products, rather than cup products.  As with cohomology determinants, the ``Massey determinants'' have slightly different constructions for closed $3$--manifolds and for $3$--manifolds with boundary.  We will give the construction for $3$--manifolds with boundary, and refer to \cite{Tur:3d} XII.2 for the construction for closed $3$--manifolds.  One can easily verify that the gluing results also apply to the Massey product determinants.

\begin{subsection}{$3$--manifolds with nonvoid boundary}
\label{subsection-massey-det-nonclosed}
To construct the determinant if $\partial M\neq\varnothing$, let $R$ be a commutative ring with
$1$, and let $K,N$ be free $R$-modules of rank $n$,$n-1$ respectively,
with $n \geq 2$ and let $S = S(K^*)$, the symmetric algebra on the
dual of $K$, as in Lemma~\ref{lemma-IntCohomDet}.  Let $f:N\times
K^{m+1} 
\rightarrow R$ be an $R$-map, with $m\geq 1$.  Define $g:N\times K
\rightarrow S$ by \[g(x,y) = \sum\limits_{i_1,\dots,i_m = 1}^n
f(x,y,a_{i_1},\dots ,a_{i_m})a_{i_1}^*\cdots a_{i_m}^* \in S\] where
$\{a_i\}_{i=1}^n$ is a basis for $K$ and $\{a_i^*\}$ is its dual
basis.  This definition for $g$ looks dependent on the basis chosen, however one can easily show independence.

Let $f_0:N\rightarrow S$ be defined by \[f_0(x) =
\sum\limits_{i_1,\dots,i_{m+1} = 1}^n f(x,a_{i_1},\dots
,a_{i_{m+1}})a_{i_1}^*\cdots a_{i_{m+1}}^* \in S.\]  Again, $f_0$ does
not depend on the chosen basis, by precisely the same argument.  Then
we have the following lemma: 
\begin{lemma}
Suppose $f_0=0$.  Let $a = \{a_i\} , b = \{b_j\}$ be bases of $K,N$
respectively, and let $\theta$ be the $(n-1\times n)$ matrix over $S$
defined by $\theta_{i,j} = g(b_i,a_j)$.  Then there exists a unique $d
= d(f,a,b) \in S^{m(n-1)-1}$ such that 
\begin{equation}
\det(\theta(i)) = (-1)^i a_i^* d.
\label{equation-IntMasseyDet}
\end{equation}
Furthermore, if $a',b'$ are other bases for $K,N$ respectively, then
\begin{equation}
d(f,a',b') = [a'/a][b'/b]d(f,a,b).
\label{equation-IntMasseyCob}
\end{equation}
\label{lemma-IntMasseyDet}
\end{lemma}
\begin{proof}
This is very similar to the proof of Lemma~\ref{lemma-IntCohomDet}.
Let $\beta$ be 
the matrix over $S$ given by $\beta_{i,j} = g(b_i,a_j)a_j^*$.  Then
the sum of the columns of $\beta$ is zero; the $i^\text{th}$ entry in
that sum is $\sum\limits_{j=1}^n \beta_{i,j} = f_0(b_i) = 0$ since our
assumption is $f_0 = 0$.  Now the same argument as given in
Lemma~\ref{lemma-IntCohomDet} to prove \eqref{equation-IntCohomDetDef}
completes the proof of 
\eqref{equation-IntMasseyDet}, and the
argument given to prove \eqref{equation-IntCohomCob} can be used to prove
\eqref{equation-IntMasseyCob}.
\end{proof}
Note that as before, over $\Z$ the determinant is well defined up to
sign, and that one may also sign-refine this determinant to remove the
sign dependence. 

We may also define the condition that $f$ is ``alternate'' in the $K$
variables; let $\overline{f_0}:N\times K \rightarrow R$ be the $R$-map
given by $\overline{f_0}(x,a) = f(x,\overbrace{a,a,\dots,a}^{m+1
\text{times}})$.  Then $f_0(x) = 0$ for all $x$ clearly implies
$\overline{f_0}(x,a) = 0$ for all $x\in N,a\in K$.  The converse is
also true provided that every polynomial over $R$ which only takes on
zero values has all zero coefficients (this is true, for example, if
$R$ is infinite with no zero-divisors).  To see why, consider $f_0(x)$
as a polynomial over $R$ ($f_0(x) \in S$ which is isomorphic to the
polynomial ring $R[a_1^* , \dots , a_n^*]$) and evaluate on the
element $(r_1 , \dots , r_n) \in R^n$; denote by $\alpha$ the
resulting element of $R$.  Then
\begin{eqnarray*}
\alpha & = & \sum_{i_1,\dots ,i_{m+1} = 1}^n f(x,a_{i_1} , \dots ,
a_{i_{m+1}}) r_{i_1}\cdots r_{i_{m+1}} \\
 & = & \sum_{i_1,\dots ,i_{m+1} = 1}^n f(x ,
r_{i_1} a_{i_1} , \dots , r_{i_{m+1}}
a_{i_{m+1}}) \\
 & = & f\left(x , \sum_{i_1= 1}^n r_{i_1}a_{i_1},\dots,\sum_{i_{m+1}=1}^n r_{i_{m+1}}a_{i_{m+1}}\right ) \\
 & = & 0.
\end{eqnarray*}
The last equality holds since all of the entries after the first are
identical.

The rest of the argument is very similar to the argument in \cite{Tur:3d}, section XII.2.  Let $M$ be a 3-manifold with nonempty
boundary, and for 
$u_1,u_2,\dots , u_k \in H^1(M)$, let $\langle u_1,\dots , u_k\rangle $ denote the
Massey product of $u_1,\dots ,u_k$ as a subset of $H^2(M)$ (note in
general this set may well be empty).  See \cite{Kra:Mas} and
\cite{Fenn:Tech} for definitions and properties of Massey products.
Now assume that $m\geq 1$ is an integer such that
\begin{center}
$(*)_m$:\quad for every $u_1,\dots,u_k\in H^1(M)$ with $k\leq m,
\langle u_1,\dots,u_k\rangle =0$
\end{center}
Here $\langle u_1,\dots,u_k\rangle =0$ means that $\langle u_1,\dots,u_k\rangle$ consists of the
single element $0\in H^2(M)$.  This condition guarantees that for any
$u_1,\dots u_{m+1}\in H^1(M)$, the set $\langle u_1,\dots,u_m\rangle$ consists of a
single element; see \cite{Fenn:Tech} Lemma 6.2.7.  Define a $\Z$-map
$f: H^1(M,\partial M) \times (H^1(M))^{m+1} \rightarrow \Z$ by 
\[ f(v , u_1 , \dots , u_{m+1}) = (-1)^m\left\langle v \cup \langle u_1,\dots,u_{m+1}\rangle  ,
[M]\right\rangle. \]
The outermost $\langle , \rangle$ is used to denote the evaluation pairing.
\begin{lemma}
$f_0 = 0$, so $f$ has a well-defined determinant (with the sign
refinement as above).
\end{lemma}
For $m=1$, condition $(*)_m$ is void, and in fact the Massey product $\langle u_1,u_2\rangle =-u_1\cup u_2$, so this reduces to Lemma~\ref{lemma-IntCohomDet}.
\begin{proof}
By the argument above, we only need to show that $f$ is alternate.
But this follows from \cite{Kra:Mas} Theorem~15, which gives that for
any element $a\in H^1(M)$, the $m+1$ times Massey product of $a$ with itself,  $\langle\overbrace{a,\dots,a}^{m+1 \text{times}}\rangle,$ lies in
$\Tors(H^2(M))$, hence cupping with an element of $H^1(M,\partial M)$ will give
an element of $\Tors(H^3(M,\partial M))$, which is null.
\end{proof}
We will call this determinant $\Det(f)$, or if we care to introduce the 
sign-refined version with a homology orientation $\omega$, $\Det_{\omega}(f)$.
Since the change of basis formula \eqref{equation-IntMasseyCob} is identical
to the change of basis formula \eqref{equation-IntCohomCob}, the sign
refinement by homology orientation is the same.
\end{subsection} 
\end{section}
%
%
\bibliographystyle{halpha}
\bibliography{bibList}
\end{document}